\documentclass{amsart}

\usepackage{amssymb}
\usepackage[all]{xy}
\usepackage{hyperref}

\usepackage{enumitem}   %from SK, fixes \enumerate

\setlist[enumerate]{itemsep=.2em,topsep=.2em,leftmargin=1.25em,itemindent=2.0em}

\newtheorem{thm}{Theorem}%[section]

\newtheorem{lem}[thm]{Lemma}
\newtheorem{cor}[thm]{Corollary}

\theoremstyle{definition}
\newtheorem{defn}[thm]{Definition}

\newtheorem{say}[thm]{}
\newtheorem{exmp}[thm]{Example}

\newtheorem{rem}[thm]{Remark}          

\newtheorem*{ack}{Acknowledgments}      % \renewcommand{\theack}{} 

\newtheorem{defn-thm}[thm]{Definition--Theorem}  %!!!!!!!!!!!!!!!!!!!!!!!!
\newtheorem{defn-lem}[thm]{Definition--Lemma}  %!!!!!!!!!!!!!!!!!!!!!!!!
\newtheorem{comments}[thm]{Comments}

\theoremstyle{remark}

\renewcommand{\c}[0]{{\mathbb C}}  

\renewcommand{\o}[0]{{\mathcal O}} 
\newcommand{\z}[0]{{\mathbb Z}}
\newcommand{\n}[0]{{\mathbb N}}
\renewcommand{\r}[0]{{\mathbb R}} 

\renewcommand{\a}[0]{{\mathbb A}}

\newcommand{\s}[0]{{\mathbb S}}

\newcommand{\p}[0]{{\mathbb P}}

\newcommand{\q}[0]{{\mathbb Q}}
\newcommand{\map}[0]{\dasharrow}
\newcommand{\qtq}[1]{\quad\mbox{#1}\quad}
\newcommand{\spec}[0]{\operatorname{Spec}}

\newcommand{\mult}[0]{\operatorname{mult}}

\newcommand{\tsum}[0]{\textstyle{\sum}}

\def\loccoh#1.#2.#3.#4.{H^{#1}_{#2}(#3,#4)}

\DeclareMathAlphabet{\mathchanc}{OT1}{pzc}%
                                {m}{it}

\newcommand{\GL}{\mathrm{GL}}

\usepackage[all]{xy}\xyoption{dvips}

\newcommand{\tprod}[0]{\textstyle{\prod}}

\begin{document}
\bibliographystyle{amsalpha}

\title[Real analytic functions]{Real analytic functions and monomial curves}
  \author{J\'anos Koll\'ar}

  \begin{abstract} We show that a function is real  analytic at the origin iff it is   arc-analytic, has a  subanalytic graph, and its  restriction to every monomial curve is analytic. This complements recent results of Kucharz and  Kurdyka.

\end{abstract}

\maketitle

%%xxx

A classical theorem of Hartogs says that a function on $\c^n$ is complex analytic iff  its restriction to every complex line parallel to  a coordinate axis is complex analytic  \cite{MR1511365}.
For functions on $\r^n$ the situation is considerably more complicated.

One can check real analyticity on 2-planes. That is, a
function  $F(x_1,\dots, x_n)$ is real analytic if its restriction to every 2-plane is real analytic  \cite{bosi71, MR4025005}; see \cite{MR4044678} for a similar result for 2-spheres.

However, restrictions to lines carry very little information,
 as shown by the  simple example   $xy/(x^2+y^2)$ on $\r^2$.  

Following \cite{MR967023}, a function  $F(x_1,\dots, x_n)$ is called {\it arc-analytic} iff the composite $F\bigl(\phi_1(t), \dots, \phi_n(t)\bigr)$  is real analytic, whenever the $\phi_i$ are real analytic. Every real analytic function is arc-analytic, but
 there are arc-analytic functions that are not even continuous; see \cite{MR1072083, MR1141439, MR1282776}.

The situation improves if we impose additional restrictions on $F$. 
\cite{MR1062969} shows that if $F$ is arc-analytic and its graph is {\it subanalytic,} then $F$ becomes  real analytic after a suitable sequence of (local, real analytic) blow-ups; see Theorem~\ref{1.4.BM.thm} for the statment and \cite{MR972342} for the basic theory of subanalytic sets. It is, however, not clear how to convert this into a
criterion for real analyticity on  $\r^n$ itself.

The recent paper \cite{ku-ku-22} proposes an  approach that looks more carefully at some special curves. Their main result says that if $F(x_1,\dots, x_n)$ is arc-analytic and
its graph  is subanalytic, then 
 $F$ is real analytic at the origin iff
$$
  F\bigl(a_1t^p+b_1t^q,\dots, a_nt^p+b_nt^q\bigr)\in \r[[t^p, t^q]]\subset
  \r[[t]]
  %%  \eqno{(\ref{reg-mon.sa.thm}.1)}
  $$
  for every $a_i, b_i\in \r$ and $0<p,q\in \n$. That is,
  $F\bigl(a_1t^p+b_1t^q,\dots, a_nt^p+b_nt^q\bigr) $ has a convergent Taylor expansion that involves only the $t$-powers  $t^{rp+sq}$ for $r,s\in\n$. 

  The aim of this paper is to prove the following variant of this result. 

\begin{thm}\label{reg-mon.sa.thm} Let $F(x_1,\dots, x_n)$ be an arc-analytic function on an open subset $U\subset \r^n$ whose graph is subanalytic. Then $F$ is analytic at the origin iff its restriction to every monomial curve is analytic  at the origin. That is, iff
  $$
  F\bigl(t^{m_1},\dots, t^{m_n}\bigr)\in \r[[t^{m_1},\dots, t^{m_n}]]\subset
  \r[[t]]
    \eqno{(\ref{reg-mon.sa.thm}.1)}
  $$
  for every $0<m_i\in \n$.
\end{thm}

Arc-analyticity implies that  $F\bigl(t^{m_1},\dots, t^{m_n}\bigr)\in
\r[[t]]$. We may always assume that $(m_1,\dots, m_n)=1$, and then
$\r[[t^{m_1},\dots, t^{m_n}]]\supset t^r\r[[t]]$ for some
$r=r(m_1,\dots, m_n)$.
(See Definition~\ref{reg.loc.rem} for explicit bounds on $r(m_1,\dots, m_n)$.)
So (\ref{reg-mon.sa.thm}.1) states the vanishing  of finitely many terms of the Taylor expansion of $F\bigl(t^{m_1},\dots, t^{m_n}\bigr)$, depending only on
$m_1,\dots, m_n$.

  \medskip

Proof.
If $F$ is analytic then locally it can be written as a power series in $x_1,\dots, x_n$, so
$F\bigl(t^{m_1},\dots, t^{m_n}\bigr)$ is a power series in $t^{m_1},\dots, t^{m_n}$, as needed.

For the converse, assume that (\ref{reg-mon.sa.thm}.1) holds.
We define elementary monomial coordinate changes in Definition~\ref{mln.coo.ch}, and
note that (\ref{reg-mon.sa.thm}.1) is preserved by them. There are now 2 steps.

(Step \ref{reg-mon.sa.thm}.2)  $F$ becomes analytic after a sequence of elementary monomial coordinate changes.  This is proved in Paragraph~\ref{main.step.pf}, using  results of \cite{hir} and \cite{MR1062969} on resolution of singularities.

Note that \cite{hir} is used through the 
Equivariant Principalization
 Theorem~\ref{hir.equiv.prin.thm}, but   only for monomial ideals, which is a relatively easy special case.
However, the proof also uses
Theorem~\ref{1.4.BM.thm}, which in turn relies on Hironaka's uniformization  theorem for subanalytic sets and his flattening theorem for proper, real analytic maps
 \cite{ MR0477121, MR0377101}; these are considerably harder. 

(Step \ref{reg-mon.sa.thm}.3)  If $F$ satisfies (\ref{reg-mon.sa.thm}.1) and becomes analytic after {\em one} elementary monomial coordinate change, then $F$ is analytic. The 2 variable case is in \cite{ku-ku-22}, the straightforward generalization is discussed  in Paragraph~\ref{one.step.say}. 

We are now done by induction on the number of elementary monomial coordinate changes. \qed

\begin{defn}[Monomial coordinate changes]\label{mln.coo.ch}
  By an {\it elementary monomial coordinate change} we mean: permuting the coordinates, choosing $1<r\leq n$ and then setting
  $\sigma:\r^n({\mathbf x'})\to \r^n({\mathbf x})$  given by
  $$
  \sigma: (x'_1, \dots,  x'_n)\mapsto
  (x'_1x'_r, \dots, x'_{r-1}x'_r, x'_r, \dots, x'_n).
  \eqno{(\ref{mln.coo.ch}.1)}
  $$
  Its inverse is defined for $(x_r\neq 0)$:
  $$
  x'_1=x_1/x_r, \dots, x'_{r-1}=x_{r-1}/x_r, x'_r=x_r, \dots, x'_n=x_n.
  \eqno{(\ref{mln.coo.ch}.2)}
  $$
  {\it Comment.} The  $x'_i$  are the usual coordinates on one of the charts after blowing up
  $(x_1=\cdots=x_r=0)$. We are guided by blow-ups, but no knowledge of them is needed if someone is willing to accept Theorem~\ref{reg-mon.bu.res.thm.0}.
  %% (\ref{reg-mon.bu.res.thm}).
  
  Thus  $F':=F\circ\sigma$ is
  $$
  F'(x'_1,\dots, x'_n):=F(x'_1x'_r,\dots, x'_{r-1}x'_r, x'_r,\dots, x'_n).
  $$
  If $\rho': t_i\mapsto (x'_i:=t_i^{m_i})$ is a monomial curve then its image
  $$
  \rho: t_i\mapsto \bigl(x_1=t^{m_1+m_r}, \dots, x_{r-1}:=t^{m_{r-1}+m_r}, x_r=t^{m_r}, \dots, x_n=t^{m_n}\bigr)
  $$
  is a monomial curve and  $F'\circ\rho'=F\circ\rho$.
  Thus if (\ref{reg-mon.sa.thm}.1) holds for $F\circ\rho$, then  also holds for   $F'\circ\rho'$  (but not conversely).

Thus the validity of (\ref{reg-mon.sa.thm}.1) for all  monomial curves is preserved by elementary monomial coordinate changes. 
\end{defn}

\begin{say}[Proof of Step~\ref{reg-mon.sa.thm}.3]
   \label{one.step.say}
  Assume that there is an  elementary coordinate change
as in (\ref{mln.coo.ch}.1)
such that
$$
F'(x'_1, \dots, x'_n):=F(x'_1x'_r, \dots, x'_{r-1}x'_r,  x'_r, \dots, x'_n)
$$
is analytic. Then $F'$  has a power series expansion involving monomials
$$
(x'_1)^{a_1}\cdots (x'_n)^{a_n}= x_1^{a_1}\cdots x_{r-1}^{a_{r-1}}x_r^{A_r}x_{r+1}^{a_{r+1}}\cdots x_n^{a_n},
$$
where $A_r(a_1,\dots, a_r)=a_r-(a_1+\cdots+a_{r-1})$. 
If $A_r(a_1,\dots, a_r)\geq 0$ for every monomial in $F'$, then
$F(x_1,\dots, x_n)$ is analytic as needed.

Otherwise there is a  monomial $M=(x'_1)^{a_1}\cdots (x'_n)^{a_n}$ in $F'$ for which
$A_r<0$.  As in Remark~\ref{unique.sol.lineq.say}, choose $(m'_1,\dots, m'_n)$ such that
$\tsum c_im'_i=\tsum a_im'_i$ has no nontrivial solution in $c_i\in \n^n$.
Set
$$
m_1:=m'_1+m'_r, \dots, m_{r-1}:=m'_{r-1}+m'_r, m_r=m'_r, \dots, m_n=m'_n.
$$
We claim that (\ref{reg-mon.sa.thm}.1) fails for this choice. Indeed, in 
$F\bigl(t^{m_1},\dots, t^{m_n}\bigr)$ we get a monomial $M\bigl(t^{m_1},\dots, t^{m_n}\bigr)$ of degree
$$
\tsum a_im'_i=-(a_1+\cdots+a_{r-1})m_r+\tsum a_im_i.
$$
Any other monomial $(x'_1)^{c_1}\cdots (x'_n)^{c_n}$ gives different $t$-degree, so $M\bigl(t^{m_1},\dots, t^{m_n}\bigr)$ is not canceled. Also,  $M\bigl(t^{m_1},\dots, t^{m_n}\bigr)\in \r[[t^{m_1},\dots, t^{m_n}]]$ iff
$$
\tsum c_im_i=-(a_1+\cdots+a_{r-1})m_r+\tsum a_im_i=\tsum a_im'_i
$$
has a solution in $c_i\in\n$. Equivalently,
$$
(c_1+\cdots+c_{r-1})m'_r+\tsum c_im'_i=\tsum a_im'_i.
\eqno{(\ref{one.step.say}.1)}
$$
On  the left hand side the coefficient of $m'_r$ is 
$c_1+\cdots+c_{r-1}+c_r$, which is at least $c_1+\cdots+c_{r-1}$.
On the right hand side $a_r<a_1+\cdots+a_{r-1}$, since we  assumed that
$A_r(a_1,\dots, a_r)=a_r-(a_1+\cdots+a_{r-1})<0$. Thus (\ref{one.step.say}.1) would be  a nontrivial solution of $\tsum c'_im'_i=\tsum a_im'_i$, a contradiction to our choice of the $(m'_1,\dots, m'_n)$. Thus
$F\bigl(t^{m_1},\dots, t^{m_n}\bigr)\not\in \r[[t^{m_1},\dots, t^{m_n}]]$.\qed
\end{say}

\begin{rem}\label{unique.sol.lineq.say}
  Fix $a_i\in \n^{>0}$ and $0<s_i\in\r$ linearly independent over $\q$.  Then $\tsum c_is_i=\tsum a_is_i$ has no nontrivial solution in $c_i\in \q$.
  Moreover, there are only finitely many $(b_i)\in \n^n$ such that
  $\bigl| \tsum b_is_i-\tsum a_is_i\bigr|<1$. Thus there is an $\epsilon>0$ such that if $|s'_i-s_i|<\epsilon$ then 
  $\tsum c_is'_i=\tsum a_is'_i$ has no nontrivial solution in $c_i\in \n$.
  Choose $s'_i$ rational and multiply through by the common denominator. We get
  $m_i\in \n^{>0}$ such that 
  $\tsum c_im_i=\tsum a_im_i$ has no nontrivial solution in $c_i\in \n$.

  We have a lot of freedom in choosing the $m_i$, for example we can arrange that all of them are primes.

\end{rem}

The technical step of the proof is the following special case of 
\cite{MR3907731}. We state it over an arbitrary field $K$, though for the current applications we need only $K=\r$.  
(Also, over finite fields
we should work with $\a^n_K:=\spec_KK[x_1,\dots, x_n ]$ instead on $K^n$.)

\begin{thm}[Monomial principalization]
  \label{reg-mon.bu.res.thm.0}
Let $K$ be a field and $h_0,\dots, h_r\in K[[x_1,\dots, x_n]]$ power series.
  Then there is a sequence of elementary monomial coordinate changes
  $\sigma_i:K^n\to K^n$ with composite
  $\bar \sigma: K^n({\mathbf z})\to K^n({\mathbf x})$
  given by $\bar \sigma^*x_i=m_i({\mathbf z})$,
  and a monomial $M({\mathbf z})$, such that the following hold.
  \begin{enumerate}
    \item The 
  $ \tilde h_i({\mathbf z}):=M({\mathbf z})^{-1}h_i\bigl(m_1({\mathbf z}),\dots,  m_n({\mathbf z})\bigr)$
  are power series, and
\item $\tilde h_j(0,\dots, 0)\neq 0$ for at least one $0\leq j\leq r$.
  \end{enumerate}
\end{thm}

While this is a special case of  \cite{MR3907731}, in Paragraph~\ref{reg-mon.bu.res.thm.0.pf} we give a proof that
uses only the statement of Hironaka's Principalization Theorem,   but very little  other algebraic geometry terminology or machinery.

\medskip

\begin{cor}[Monomial regularization]
  \label{reg-mon.bu.res.thm} Let $h_i$ be real analytic functions on an open set  ${\mathbf 0}\in U({\mathbf x})\subset \r^n({\mathbf x})$.
  Then there is a sequence of elementary monomial coordinate changes
  $\bar\sigma=\sigma_1\circ\cdots\circ \sigma_s:\r^n({\mathbf z})\to \r^n({\mathbf x}) $ and an open set  ${\mathbf 0}\in U({\mathbf z})\subset \r^n({\mathbf z})$,
  such that
  $$
  (h_0{:}\cdots{:}h_r)\circ \bar\sigma:
  U({\mathbf z})\map \p^r
  $$
  is  real analytic at the origin.
\end{cor}

  Note that this is different from Hironaka-type resolution of indeterminacies in 2 ways.
 It is easier, since we want the composite to be real analytic at 1 point only, not everywhere. On the other hand,   we  are allowed to use only elementary monomial coordinate changes, while
the usual resolution methods also use local analytic coordinate changes between
successive  blow-ups in an essential way.
\medskip

Proof. We apply Theorem~\ref{reg-mon.bu.res.thm.0} to the power series expansions of the $h_i$ to get $\bar \sigma: \r^n({\mathbf z})\to \r^n({\mathbf x})$
and 
 $\tilde h_i:=(h_i\circ \bar \sigma)/M$   as in (\ref{reg-mon.bu.res.thm.0}.1).  The  $\tilde h_i$ are analytic at the origin, and $\tilde h_j$ is nonzero at the origin for some $j$. Thus 
$$
  (h_0{:}\cdots{:}h_r)\circ \bar\sigma:
  U({\mathbf z})\map \p^r
  $$
  agrees with
  $$
  \bigl(\tilde h_0{:}\cdots{:}\tilde h_r\bigr):
  U({\mathbf z})\map \p^r,
  $$
and the latter is defined and analytic at the origin.\qed

\medskip

The following is proved, though differently stated, in \cite[Thm.1.4]{MR1062969}.

\begin{thm}  \label{1.4.BM.thm}
  %%\cite[Thm.1.4]{MR1062969}
  Let $X$ be a real analytic manifold and $f:X\to \r$  an  arc analytic function  whose graph is subanalytic.
    Then there is a sequence of real analytic morphisms
$$
    X_m\stackrel{\pi_m}{\longrightarrow} \tilde X_{m-1}\stackrel{u_{m-1}}{\longrightarrow}X_{m-1}\stackrel{\pi_{m-1}}{\longrightarrow}\cdots \stackrel{\pi_1}{\longrightarrow}\tilde X_0\stackrel{u_1}{\longrightarrow} X_0=X,
    \eqno{(\ref{1.4.BM.thm}.1)}
  $$
  where, for every $i$,
  \begin{enumerate}\setcounter{enumi}{1}
  \item  $u_i: \tilde X_i\to X_i$ is a locally finite, open  covering, 
\item $\pi_i: X_i\to \tilde X_{i-1}$ is the blow-up of a smooth, closed, real analytic subset, and
  \item the composite
  $F_m:=F\circ u_1\circ \pi_1\circ\cdots\circ u_{m-1}\circ \pi_m: X_m\to \r
  $
    is real analytic.  \qed
    \end{enumerate}
\end{thm}

\begin{say}[Proof of Step \ref{reg-mon.sa.thm}.2] \label{main.step.pf}
  Let $F:U_0\to \r$ be as in Theorem~\ref{reg-mon.sa.thm}, and   consider the sequence
  of morphisms (\ref{1.4.BM.thm}.1).

  We prove that there is a
  sequence of elementary monomial coordinate changes
  $\sigma_i:\r^n\to \r^n$ and an open neighborhood $0\in U_m\subset \r^n$, such that   
  $$
  \tau_0:=\sigma_1\circ\cdots\circ \sigma_{s(m)}:
  U_m\to \r^n=:X_0
  $$
  lifts to  $\tau_m:U_m\to X_m$ that is real analytic at the origin. If this holds then
  $ F\circ \tau_0=F_m\circ \tau_m$ is real analytic at the origin, as needed.

  We construct the lifting $\tau_m$ by going up the chain (\ref{1.4.BM.thm}.1).
    Thus assume that we already have $\tau_i:U_i\to X_i$.
  Since  $u_i: \tilde X_i\to X_i$ is a locally finite, open  covering, any
  $\tau_i:U_i\to X_i$ lifts to $\tilde\tau_i:\tilde U_i\to \tilde X_i$
  over some $0\in \tilde U_i\subset U_i$.

  It remains to lift each $\tilde\tau_i:\tilde U_i\to \tilde X_i$
  to $\tau_{i+1}:U_{i+1}\to X_{i+1}$.

  Note that, in a neighborhood of  $\tilde\tau_i({\mathbf 0})$,  the blow-up $\pi_{i+1}$ is the closed graph of the meromorphic map
  $$
  (h_0^{(i)}{:}\cdots{:}h_{r_i}^{(i)}): \tilde X_{i}\map \p^{r_i},
  $$
  where the    $h_j^{(i)}$ are real analytic.
  Consider the composite
  $$
  H^{(i)}:=\bigl(h_0^{(i)}\circ \tilde\tau_i{:}\cdots{:}h_{r_i}^{(i)}\circ \tilde\tau_i\bigr): \tilde U_{i}\map \p^{r_i}.\eqno{(\ref{main.step.pf}.1)}
  $$
  By Corollary~\ref{reg-mon.bu.res.thm} there is a   sequence of elementary monomial coordinate changes
  $\sigma^{(i)}:\r^n\to \r^n$  and $U_{i+1} \subset \r^n$, such that
  $H^{(i)}\circ \sigma^{(i)}: U_{i+1}\to \p^{r_i}$ is
  real analytic at the origin. Then
  $$
  \tau_{i+1}:=\bigl(\tilde\tau_i\circ \sigma^{(i)}, H^{(i)}\circ \sigma^{(i)}\bigr): U_{i+1}\to X_{i+1}\subset \tilde X_i\times \p^{r_i}
  $$
  is real analytic at the origin.
   \qed

  \end{say}

\subsection*{Further remarks and examples}{\ }

\begin{rem} \label{bezout.est.say}
  Even over $\c$, one can not check analyticity at a point using bounded degree curves only.
  
To see this, let $h(x,y)$ be any irreducible polynomial of degree $n$. Set 
  $F=H(x,y)/h(x,y)$ where $H$ is any analytic function such that
  $\mult_{(0,0)}H>n(n-1)$.  As we discuss below,  the restriction of $F$ to every curve of degree $<n$  through the origin is analytic  at the origin, but $F$ is usually not.

To see this let $C:=(g=0)\subset \c^2$ be a curve of degree  $d<n$. Then $h$ does not vanish on $C$. Set $R:=\c[x,y]/(g)$. Multiplication by $h$ gives an isomorphism
  $$
  R[h^{-1}]/R \cong R/(h) \cong \c[x,y]/(g,h),
  $$
  and the latter has dimension $\leq nd$ by B\'ezout's theorem.
  Both $x, y$ are nilpotent on the summand supported at the origin, so
  $x^iy^j/h$ is regular on $C$ at the origin whenever $i+j>nd$.

  It would also not be enough to use bounded multiplicity (but unbounded degree) monomial curves in
Theorem~\ref{reg-mon.sa.thm}. For example, set $F_m=x^my^m/(x^2+y^2)$.
Using Definition~\ref{reg.loc.rem} we see that
$$
F_m(t^{m_1}, t^{m_2})\in \c[[t^{m_1}, t^{m_2}]]
\qtq{if} \min\{m_1,m_2\}\leq m/2.
  $$ 
On the other hand, it is clear from the proof of Theorem~\ref{reg-mon.sa.thm}, especially Paragraphs~\ref{one.step.say}--\ref{unique.sol.lineq.say},  that many sparse subsets of all monomial curves do work.
For example, as noted in Remark~\ref{unique.sol.lineq.say}, it is enough to check  (\ref{reg-mon.sa.thm}.1) when the $m_i$ are all primes. 
\end{rem}

\begin{defn}[Frobenius number]\label{reg.loc.rem}
The largest $m$ such that   $t^m$ is not in 
$ K[t^{m_1},\dots, t^{m_n}]$ is called the {\it Frobenius number} of
$\{m_1,\dots, m_n\}$.

For $n=2$ the Frobenius number  is  $m_1m_2-m_1-m_2$. For $n\geq 3$
even good estimates are hard, but
it is known that  the  Frobenius number  is at most $(2/n)m_1m_n$
(assuming $m_1\leq \cdots\leq  m_n$), see \cite{MR311565}.

Thus if  $F=g/h$ for some explicit power series $g,h$,
then the condition (\ref{reg-mon.sa.thm}.1) is algorithmically decidable.
\end{defn}

A purely algebraic version of Theorem~\ref{reg-mon.sa.thm} is the following.

\begin{thm}\label{reg-mon.thm0} Let $K$ be an arbitrary field
  and $f, g_i\in  K[[x_1,\dots, x_n]]$ power series. Then $f\in (g_1, \dots, g_r)$ iff
  $$
  f(t^{m_1},\dots, t^{m_n})\in
  \bigl(g_1(t^{m_1},\dots, t^{m_n}),  \dots, g_r(t^{m_1},\dots, t^{m_n})\bigr)
  \eqno{(\ref{reg-mon.thm0}.1)}
  $$
in the ring $K[[t^{m_1},\dots, t^{m_n}]]$
   for every $0<m_i\in \n$. 
\end{thm}

Note that we need (\ref{reg-mon.thm0}.1) in
the ring $K[[t^{m_1},\dots, t^{m_n}]]$, not in $K[[t]]$.
This is a somewhat subtle condition, since
$K[[t^{m_1},\dots, t^{m_n}]]$ is not a principal ideal domain.

Requiring  (\ref{reg-mon.thm0}.1) in $K[[t]]$ would lead to the notion of integral dependence, see \cite{laz-book, swa-hun, MR2499856}.

This works over any field since principalization of monomial ideals is a combinatorial question,  hence independent of the characteristic, see \cite{MR2219845}.
However, Theorem~\ref{reg-mon.thm0}
 has a simpler,  more direct proof using 
 Gr\"obner bases.

\subsection*{Proof of Theorem~\ref{reg-mon.bu.res.thm.0}}{\ }

We need the following special case of resolution of singularities. We discuss its proof at the end.
  
\begin{thm}[Monomial principalization of monomial ideals]
  \label{reg-mon.bu.res.thm.mon}
  Let $K$ be a field and
$J\subset K[x_1,\dots, x_n]$  a monomial ideal.
    Then there is a sequence of elementary monomial coordinate changes
  $\sigma_i:K^n\to K^n$ with composite
  $\bar \sigma: K^n({\mathbf z})\to K^n({\mathbf x})$
    given by $\bar \sigma^*x_i=m_i({\mathbf z})$,
    such that $\bar\sigma^*J\subset K[z_1,\dots, z_n]$ is a
    principal ideal, generated by a  monomial.
\end{thm}

\begin{say}[Proof of Theorem~\ref{reg-mon.bu.res.thm.0}]
  \label{reg-mon.bu.res.thm.0.pf}
Let $E(h_0,\dots, h_r)\subset \n^n$ be the set of exponents
$(a_1,\dots, a_n)$ for which $x_1^{a_1}\cdots x_n^{a_n}$
appears in some $h_i$ with nonzero coefficient.

The observation of  \cite{MR3907731} is that one should work with the monomial ideal
$$
J=J(h_0,\dots, h_r):= \bigl(x_1^{a_1}\cdots x_n^{a_n} : (a_1,\dots, a_n)\in E(h_0,\dots, h_r) \bigr)\subset K[x_1, \dots, x_n].
$$
Next we apply Theorem~\ref{reg-mon.bu.res.thm.mon}
to get a composite of a sequence of elementary monomial coordinate changes
$\bar\sigma:K^n\to K^n$ such that  $\bar\sigma^*J$
is a principal ideal generated by 
$$
M({\mathbf z})=m_1^{a_1}\cdots m_n^{a_n} \qtq{for some}
(a_1,\dots, a_n)\in E(h_0,\dots, h_r).
$$
Also, if  $x_1^{a_1}\cdots x_n^{a_n}$
appears in some $h_j$ with nonzero coefficient,
then  $M({\mathbf z})$  appears in the same  $\bar\sigma^*h_j$ with nonzero coefficient.

Since $M({\mathbf z})$  divides
each monomial in $\bar\sigma^*h_i$,  the 
$$
\tilde h_i({\mathbf z}):=M({\mathbf z})^{-1}h_i\bigl(m_1({\mathbf z}),\dots,  m_n({\mathbf z})\bigr)
$$
are power series,  and
$\tilde h_j(0,\dots, 0)\neq 0$. \qed
\end{say}

Next we discuss how to get Theorem~\ref{reg-mon.bu.res.thm.mon} from the
large literature on resolution of singularities.
First we outline how the method of Hironaka gives this, though here we need to
`follow the proof.'
Then we show how to derive 
Theorem~\ref{reg-mon.bu.res.thm.mon} from the statement of
the Equivariant Principalization Theorem.

\begin{say}[Hironaka's method] As in \cite{hir}, we start with an ideal
  $J:=(f_1,\dots, f_r)\subset K[x_1,\dots, x_n]$. We aim to produce a smooth subvariety $Z\subset K^n$ whose blow-up `simplifies' the singularities of $J$.
  The procedure is rather complicated, going through a sequence of new ideals with new generators. However, all new generators are constructed by repeatedly applying 2 steps:
  \begin{itemize}
  \item taking a partial derivative of a generator, or
  \item multiplying 2 generators.
  \end{itemize}
  Both of these steps preserve monomials, so 
  if we start with monomials $f_i$, then all ideals in the Hironaka process are generated by  monomials.

  At the end we arrive at an ideal $J^*:=(f^*_1,\dots, f^*_s)\subset K[x_1,\dots, x_n]$ such that $J^*$ is the ideal of a  smooth subvariety $Z\subset K^n$.

  The zero set of any monomial ideal is a union of intersections of coordinate hyperplanes, and the only smooth ones are intersections of coordinate hyperplanes.
  Thus $J^*$ defines an intersection of coordinate hyperplanes.
  The corresponding blow-ups are exactly the
   elementary monomial coordinate changes, as in
   Definition~\ref{mln.coo.ch}.

   It turns out that we are in a much simpler situation than the one considered by Hironaka, where it is quite difficult to establish that the process is independent of the choice of the generators. In our case we have natural generators.

   Also, as shown in \cite{MR2219845}, 
   the original method of Hironaka applies to
   monomial ideals in any characteristic.
   The key change is that, instead of differentiation   $(cx^m)'=mcx^{m-1}$, one should use the operator
   $cx^m\mapsto cx^{m-1}$. Then we get a process that works  in positive characteristic as well. \qed
   \end{say}

Monomial ideals are usually dealt with as special cases of toric geometry.
For us the simplest  explicit reference appears to be  the following.

\begin{thm}[Equivariant Principalization]
  \label{hir.equiv.prin.thm}
    Let $K$ be a field of characteristic 0.
    Let $G \curvearrowright X$ be a group acting on a $K$-variety and
    $I\subset \o_X$ a $G$-invariant ideal sheaf. Then there is a
    sequence
  $$
    (G\curvearrowright X_s)\stackrel{\pi_{s-1}}{\longrightarrow}
    (G\curvearrowright X_{s-1})\to \cdots \to
    (G\curvearrowright X_1)\stackrel{\pi_{0}}{\longrightarrow}
    (G\curvearrowright X_0=X)
    $$
    of blow-ups of  $G$-invariant, smooth subvarieties, such that
    $$
    I_s:=(\pi_0\circ\cdots\circ \pi_s)^*I\subset \o_{X_s}
    $$
    is locally principal. \qed
\end{thm}

{\it References.}  This is a special case of \cite[3.26]{k-res}.
There is no group action in the statement, but 
 \cite[3.26.2]{k-res}
says that the principalization given there is `functorial on smooth morphisms.'
Isomorphisms are smooth morphisms, so the principalization in \cite[3.26]{k-res} is invariant under all isomorphisms.

The theorem is slightly more general than \cite[p.142]{hir}; the improvements in  \cite{bie-mil-inv, MR1654779} imply the general case.

 \medskip

We aim to apply Theorem~\ref{hir.equiv.prin.thm} to  $X=\a^n_K$ and
$I\subset K[x_1,\dots, x_n]$ an ideal generated by monomials.
We can choose the group to be $(K^\times)^n$.

\begin{defn}   A {\it standard toric} chart over $K$ is
  $(T^n\curvearrowright K^n)$, where we view $K^n$ and $T^n:=(K^\times)^n$
  as 
  $K$-varieties, and the  action is  given by
  $$
(\lambda_1, \dots, \lambda_n)\times (x_1,\dots,x_n)\mapsto   (\lambda_1x_1, \dots, \lambda_nx_n).
$$
  We allow permuting the $x_i$ and changing the action by any $(m_{ij})\in \GL_n(\z)$ to
  $$
(\lambda_1, \dots, \lambda_n)\times (x_1,\dots,x_n)\mapsto   \bigl(\tprod_j \lambda_j^{m_{1j}}x_1, \dots, \tprod_j \lambda_j^{m_{nj}}x_n\bigr).
  $$
  A {\it toric} $K$-variety is a $(T^n\curvearrowright X)$, where $X$ is a smooth  $K$-variety of dimension $n$,  that is covered by $T^n$-invariant, standard charts.

  {\it Comment on terminology.} If $K$ is algebraically closed, this is the usual notion of smooth toric variety. If $K$ is not algebraically closed, many authors allow replacing $T^n$ with some  other $K$-form of it. Thus, over $\r$, one could use $(\s^1)^n$ instead. For us it is better to work with
  $T^n=(K^\times)^n$.
  
\end{defn}

The prime example, and the only one we need, is the following.

\begin{exmp} Let $Z\subset K^n$ be the intersection of some coordinate hyperplanes, say $(x_1=\cdots=x_r=0)$. The blow-up $B_ZK^n$ is
   toric. It is covered by $r$ standard toric charts, each one an
  elementary monomial coordinate change.

  The $r$th chart is described in Definition~\ref{mln.coo.ch}. We get the $i$th chart by permuting the $1,\dots, r$ such that $i$ ends up in the last position, and then
  using (\ref{mln.coo.ch}.1--2).

  Note that the only $T^n$-invariant subvarieties on $(T^n\curvearrowright K^n)$ are the intersections of coordinate hyperplanes.
  Therefore,  if 
  $(T^n\curvearrowright X)$ is a  toric $K$-variety, and if 
  $Z\subset X$ is a  $T^n$-invariant subvariety, then
  $(T^n\curvearrowright B_ZX)$ is also a  toric $K$-variety.
  Iterating these gives the following.
\end{exmp}

  \begin{lem} \label{lem.18} Let
    $$
    (T^n\curvearrowright X_s)\stackrel{\pi_{s-1}}{\longrightarrow}
    (T^n\curvearrowright X_{s-1})\to \cdots \to
    (T^n\curvearrowright X_1)\stackrel{\pi_{0}}{\longrightarrow}
    (T^n\curvearrowright K^n)
    $$
    be a sequence of blow-ups of  $T^n$-invariant subvarieties.
    Then $X_s$ is covered by standard toric charts $X_s^j$, such that
    each $K^n\cong X_s^j\to K^n$ is a composite of $\leq s$ 
    elementary monomial coordinate changes. \qed
  \end{lem}

  \begin{say}[Proof of Theorem~\ref{reg-mon.bu.res.thm.mon}]
  Combining Theorem~\ref{hir.equiv.prin.thm} and
  Lemma~\ref{lem.18}  gives
   a sequence of elementary monomial coordinate changes
  $\sigma_i:K^n\to K^n$ with composite
  $\bar \sigma: K^n({\mathbf z})\to K^n({\mathbf x})$
    given by $\bar \sigma^*x_i=m_i({\mathbf z})$,
    such that $\bar\sigma^*J\subset K[z_1,\dots, z_n]$ is a
    principal ideal.

    The pull-back of a monomial by a sequence of elementary monomial coordinate changes is a monomial, so  $\bar\sigma^*J\subset K[z_1,\dots, z_n]$ is a monomial ideal. Since every divisor of a monomial is a monomial, if a
 monomial ideal $\bar\sigma^*J=(M_1,\dots, M_s)$ is principal, then
its is generated by one of the $M_i$. \qed
\end{say}

 \begin{comments}[Other principalization methods]\label{comm.on.refs.say}
Principalization of monomial ideals is a combinatorial question, which can be done in many different ways; see
\cite{kkms, amrt, MR1748623, MR1892938, MR2343384}.

In our proof we rely on Step~\ref{reg-mon.sa.thm}.3, thus 
 we need a process that uses only 
 elementary blow-up sequences.  The  above  papers are not designed to do this.

 However, it would be possible to establish a version of
 Step~\ref{reg-mon.sa.thm}.3 for an arbitrary monomial coordinate change.
 Once that is done, one could use any toric principalization method instead of
 Theorem~\ref{hir.equiv.prin.thm}.
 \end{comments}

\begin{ack}  I thank  W.~Kucharz and K.~Kurdyka for introducing me to this subject, and many helpful comments and corrections.
Partial  financial support    was provided  by  the NSF under grant number
DMS-1901855.
\end{ack}

%% \bibliography{../refs-main/refs}

\def\cprime{$'$} \def\cprime{$'$} \def\cprime{$'$} \def\cprime{$'$}
  \def\cprime{$'$} \def\dbar{\leavevmode\hbox to 0pt{\hskip.2ex
  \accent"16\hss}d} \def\cprime{$'$} \def\cprime{$'$}
  \def\polhk#1{\setbox0=\hbox{#1}{\ooalign{\hidewidth
  \lower1.5ex\hbox{`}\hidewidth\crcr\unhbox0}}} \def\cprime{$'$}
  \def\cprime{$'$} \def\cprime{$'$} \def\cprime{$'$}
  \def\polhk#1{\setbox0=\hbox{#1}{\ooalign{\hidewidth
  \lower1.5ex\hbox{`}\hidewidth\crcr\unhbox0}}} \def\cdprime{$''$}
  \def\cprime{$'$} \def\cprime{$'$} \def\cprime{$'$} \def\cprime{$'$}
\providecommand{\bysame}{\leavevmode\hbox to3em{\hrulefill}\thinspace}
\providecommand{\MR}{\relax\ifhmode\unskip\space\fi MR }
\providecommand{\MRhref}[2]{%
  \href{http://www.ams.org/mathscinet-getitem?mr=#1}{#2}
}
\providecommand{\href}[2]{#2}

 \bigskip

  Princeton University, Princeton NJ 08544-1000, \

  \email{kollar@math.princeton.edu}

\end{document}